\documentclass[12pt]{amsart}
\usepackage[english]{babel}
\usepackage{amsmath}
\usepackage{amsthm}
\usepackage{amssymb}
\usepackage{mathrsfs}
\usepackage{enumerate}
\usepackage[notcite, final, notref]{showkeys}
\usepackage{amsfonts}
\usepackage{dsfont}
\usepackage{tikz}
\usepackage[siunitx]{circuitikz}
\usetikzlibrary{arrows,calc,chains,shapes,dsp}
\def\C{\mathbb C}
\def\R{{\mathbb R}}
\newtheorem{Pa}{Paper}[section]
\newtheorem{Tm}[Pa]{{\bf Theorem}}
\newtheorem{La}[Pa]{{\bf Lemma}}
\newtheorem{Ob}[Pa]{{\bf Observation}}
\newtheorem{Cy}[Pa]{{\bf Corollary}}
\newtheorem{Rk}[Pa]{{\bf Remark}}
\newtheorem{Pn}[Pa]{{\bf Proposition}}

\newtheorem{Ex}[Pa]{{\bf Example}}
\newtheorem{Dn}[Pa]{{\bf Definition}}
\usepackage{xcolor}

\DeclareUnicodeCharacter{2212}{-}

\title[MCIC \& Passive Systems]
{Passive Linear Continuous-time Systems:\\[0.2cm]
Characterization through Structure}
\author[I. Lewkowicz]{Izchak Lewkowicz}
\address{(IL) School of Electrical and Computer Engineering 
Ben-Gurion University of the Negev\\ P.O.B. 653\\ Beer-Sheva, 84105\\
Israel}
\email{izchak@bgu.ac.il}

\begin{document}
\bibliographystyle{plain}
\begin{abstract}
We here show that the family of finite-dimensional, continuous-time,
passive, linear, time-invariant systems can be characterized through
the structure of maximal matrix-convex cones, closed under inversion.
Moreover, this observation unifies three setups:\\
(i) differential inclusions,\\
(ii) matrix-valued rational functions,\\
(iii) realization arrays associated with rational functions.
\vskip 0.2cm

\noindent
It turns out that in the discrete-time case, the corresponding
structure is of a maximal matrix-convex set, closed under
multiplication among its elements.
\vskip 0.2cm
\end{abstract}
\maketitle


\date{today}
\tableofcontents

\bibliographystyle{plain}
\section{Introduction}
\setcounter{equation}{0}

\noindent
In the study of dynamical systems, passivity is a fundamental property.
Thus, it has been extensively studied in various frameworks. For a
modest sample of the vast literature on the subject, see e.g.
\cite{AnderMoore1968}
\cite{AnderVongpa1973},
\cite{Belev1968},
\cite{BottDuff1949}-\cite{Cauer1932},
\cite{CohenLew2007},
\cite{Duffin},
\cite{Lewk2020c},
\cite{Lewk2020b},
\cite{MorelliSmith2019},
\cite{Smith2002}
\cite{Willems1st1972}-\cite{Willems1976}
and for the infinite-dimensional case, see e.g. \cite{BallStaff2006},
\cite{Staff2002} and \cite{Wohl1969}.
\vskip 0.2cm

\noindent
Here we confine the discussion to finite-dimensional, linear time-invariant
and continuous-time. Then, passive systems are modelled by {\em Positive
Real} functions, denoted by $~\mathcal{PR}$, namely (in scalar setup)
rational functions which analytically map $\C_R$, the open right half plane,
to its closure $\overline{\C}_R$. Namely the family of real scalar rational
function $f(s)$ of a complex variable $s$, where ${\rm Re}(f)\geq 0$,
whenever ${\rm Re}(s)>0$. 
\vskip 0.2cm

\noindent
For example, a scalar rational function of McMillan degree one, is
positive real, if and only if it is of the form of
\begin{equation}\label{eq:Deg1}
{\rm either}\quad f(s)=a+bs\quad{\rm or}\quad
h(s)=d+\frac{b}{s+a}\quad\quad
b>0,~a, d\geq 0.
\end{equation}
($f(s)$, $h(s)$ are $\mathcal{PR}$ functions  of degree zero, when $b=0$).
\vskip 0.2cm

\noindent
The fact that the set $~\mathcal{PR}~$ may serve as a prototype model to
continuous-time, linear, passive systems is classical, see e.g.
\cite[Theorem 2.6.1]{AnderVongpa1973}, \cite[Section 3.18]{Belev1968},
\cite{MorelliSmith2019}, \cite{Smith2002} and \cite{Willems2nd1972}.
\vskip 0.2cm

\noindent
Probably the better known physical realization of $~\mathcal{PR}$ functions
is through electrical circuits comprised of resistors, inductors and
capacitors ($R-L-C$ circuits).
Specifically, recall that in 1931 O. Brune, \cite{Brune1} showed the
following:\\
{\em 
The driving point immittance of a lumped\begin{footnote}{Impedance of
{\em distributed} $R-L-C$ networks are associated with {\em irrational}
$~\mathcal{PR}$ functions, see e.g. \cite{Wohl1969}.}\end{footnote}
\mbox{$R-L-C$} electrical network is a rational $~\mathcal{PR}$
function.}
\vskip 0.2cm

\noindent
In 1949 R. Bott and R.J. Duffin, \cite{BottDuff1949}, showed that also
the converse is true:\\
{\em An arbitrary positive real rational function can be realized as the
driving point immittance of a lumped \mbox{$R-L-C$} electrical network.
}
\vskip 0.2cm

\noindent
For example, the rational function $h(s)$ in Eq. \eqref{eq:Deg1} can be
realized as the driving point impedance of the simple circuit in Figure
\ref{Fig:CisDrivingPointImpedance2}

\begin{figure}[ht!]
\begin{tikzpicture}[scale=1.8]
  \draw[color=black, thick]
     (0,0) to [short,o-] (3,0){} % Baseline for connection to ground
     (-0.1,0.8) node[]{\large{$\mathbf{Z_{\rm in}~~\rightarrow}$}}
     (0,1.6) to [short,o-] (0.1,1.6)
     (0.1,1.6)  to [R,l=$R_1$,](2,1.6)
     (2,1.6)   to node[short]{} (3,1.6)
     (3,0) to [C,l=$C$, *-*] (3,1.6)
     (2,0) to [R, l=$R_2$, *-*] (2,1.6)
     ;
\end{tikzpicture}
\caption{
$Z_{\rm in}=d+\frac{b}{s+a}
\quad {\scriptstyle R_1}=d~~
{\scriptstyle\frac{1}{C}}=b~~
{\scriptstyle\frac{1}{R_2C}}=a
$
}
\label{Fig:CisDrivingPointImpedance2}
\end{figure}
\vskip 0.2cm

\noindent
Duality between rational positive real functions and the driving point
immittance of $~R-L-C~$ electrical circuits, has already been long
recognized, e.g. \cite{Brune1}, \cite{Cauer1926}, \cite{Cauer1932}.
This has lead to rich and well-established theory, see e.g.
\cite{AnderVongpa1973}, \cite{BottDuff1949}, \cite{Belev1968},
\cite{Duffin}, \cite{Goldberg1962} and \cite{Wohl1969}. For a recent
comprehensive account of circuits describing
$~\mathcal{PR}$ functions of degree two, see \cite{MorelliSmith2019}.
\vskip 0.2cm

\noindent
Already at this stage, note that the family $\mathcal{PR}$ is closed
under the following operations: (i) positive scaling, (ii) summation
and (iii) inversion, see e.g. \cite[Theorem 1.1]{Goldberg1962}.\\
In fact these operations have an electrical circuits
interpretation:
\begin{center}
$\begin{matrix}
{\rm
Positive~scaling}&-&{\rm
transformer~ratio}\\
{\rm
Summation}&-&{\rm
series~connection~of~impedances}\\
{\rm
Inversion}&-&{\rm
impedance / admittance}.
\end{matrix}$
\end{center}
\vskip 0.2cm

\noindent
An alternative physical realization of $~\mathcal{PR}$ functions is
given by the classical analogy between \mbox{$R-L-C$} electrical
circuits and simple mechanical systems, see e.g.\begin{footnote}
{where it is explained why ``inerter" replaces
 ``mass".}\end{footnote} \cite{Smith2002}.

\begin{center}
$\begin{array}{c|c}
{\rm
electrical}~~&~~{\rm
mechanical}\\
\hline\\
{\rm  current}~~&~~{\rm force}\\~
{\rm voltage}~~&~~{\rm velocity}\\
{\rm transformer}~~&~~{\rm gear~~transmission}\\
{\rm resistor~(admittance)}~~&~~{\rm damper}\\
{\rm inductor~(admittance)}~~~&~~{\rm spring}\\
{\rm capacitor~(admittance)}~~&~~{\rm inerter}.
\end{array}$
\end{center}
\vskip 0.2cm

\noindent
Following the above-mentioned intuition that the impedance of
all single-input single-output $R-L-C$ circuits forms a convex
cone, closed under inversion, in a series of papers
\cite{AlpayLew2011}, \cite{AlpayLew2013}, \cite{AlpayLew2018},
\cite{AlpayLew2019a}, \cite{CohenLew1997a}, \cite{CohenLew1997b},
\cite{CohenLewRod1997}, \cite{Lew1999a} and \cite{LewRodYar2005},
we have explored the connection between this structure and the family
of finite-dimensional, continuous-time passive, linear time-invariant
scalar systems. In the sequel we show that the passage to multi-input
multi-output $R-L-C$ circuits, lead to the substitution of the notion
of classical convexity by (the more restrictive) {\em matrix-convexity}.
\vskip 0.2cm

\noindent
There has been characterizations of passive systems in various setups.
As already mentioned, the list of references here is only a modest
sample of the vast literature on the subject. Probably the most intuitive
description of passivity is based on the notion of ``storage function",
due to J.C. Willems, see e.g. \cite{Willems1st1972}.
\vskip 0.2cm

\noindent
Here we adopt a more abstract point of view and focus on the following
question:
\begin{quote}
How can one characterize the family of finite-dimensional,
continuous-time, passive, linear time-invariant systems 
%(of prescribed dimensions) 
through the structure of the whole set.
\end{quote}
%\vskip 0.2cm

\noindent
The answer is that this family forms a maximal {\em matrix-convex}
cone, closed under inversion. Moreover, this observation unifies
three setups:
\begin{itemize}

\item[(i)~~]{}Differential inclusions,

\item[(ii)~]{}Positive real rational functions,

\item[(iii)]{}Families of realization arrays of positive real
rational functions.
\end{itemize}
\vskip 0.2cm

\noindent
This work is organized as follows. In Section \ref{sec:DiffInc}
we present the first motivation: Stability of differential
inclusions. The background foundations of {\em matrix-convex} sets
and cones, to be used in the sequel, are laid in Section
\ref{Sec:Mcic}. Then this concept is applied to maximal non-singular
matrix-convex cones, closed under inversion: of (constant) matrices
and of matrix-valued rational functions, in Sections
\ref{Sec:Maximal Mcic} and \ref{MCICs of Rational Functions},
respectively. At the second part of Section \ref{MCICs of Rational
Functions} an application of the first part to the design of
elaborate feedback-loop systems, is introduced. Finally in Section
\ref{Sec:Realizations} we explore realization arrays of families
of systems.

\section{stability of differential inclusions}
\label{sec:DiffInc}
\setcounter{equation}{0}

\noindent
As a first motivation we resort to the problem of stability of
differential inclusions, see e.g. \cite{CohenLewRod1997},
\cite{MolcPyat1989}. In engineers circles it is referred to as
``robust exponential stability" (and colloquially as
``quadratic stability").
\vskip 0.2cm

\noindent
For a set \mbox{$\mathbf{A}=\left\{A_1,~\ldots~,~ A_m\right\}$}
of $n\times n$ matrices let the differential inclusion,
\begin{equation}\label{eq:DiffInclusion}
\frac{dx}{dt}\in\mathbf{A}x
\quad\quad\quad x\in\mathbb{R}^n,
\end{equation}
mean that there exists an ~{\em unknown}~ selection
\mbox{$A(t)\in\mathbf{A}$}, \mbox{$t\geq t_o$} and initial value
\mbox{$x(t_o)=x_o$}, so that $x(t)$ is a solution of the system
\begin{equation}\label{eq:selection}
\frac{dx}{dt}=A(t)x(t)
\quad\quad x(t_o)=x_o~,~~t_o\geq 0.
\end{equation}
Although the precise nature of $A(t)$ is not necessary for the
sequel, yet for the sake of precision we adopt the following
standard assumption: (i) All selections $A(t)$ are Lebesgue
measurable and locally integrable. Thus, piecewise constant
and (other discontinuities) are allowed. (ii) All solutions
$x(t)$ are absolutely continuous. (iii) The equality
\mbox{$\frac{dx}{dt}=A(t)x(t)$} holds almost
everywhere\begin{footnote}
{Under this assumptions, for a particular selection,
Eq. \eqref{eq:selection} is equivalent to the integral equation
$x(t)=x(t_o)+\int_{t_o}^tA(\tau)x(\tau)d(\tau)$.}\end{footnote}.
\vskip 0.2cm

\noindent
To state conditions for stability of the above differential inclusion
we need to introduce some notations:
$~\overline{\mathbf H}_n$  will denote the set of
(possibly singular) $n\times n$ Hermitian matrices and
$~{\mathbf H}_n$  will be the subset of non-singular Hermitian matrices.
The respective
subsets of positive (semi)-definite matrices are denoted by
($\overline{\mathbf P}_n$) and ${\mathbf P}_n$.
Let also $i\overline{\mathbf H}_n~$ be the set of
\mbox{$n\times n$} Skew-Hermitian matrices.
\vskip 0.2cm

\noindent
Next, consider the set of $n\times n$ matrices $A$ all
satisfying a Lyapunov inclusion with the same factor,
\begin{equation}\label{eq:DefL_H}
\begin{matrix}
\mathbf{L}_H&:=&\{ A~:~HA+A^*H\in\mathbf{P}_n~\}
\\~\\
\overline{\mathbf L}_H&:=&\{ A~:~HA+A^*H\in\overline{\mathbf P}_n~\}
\end{matrix}
\quad\quad{\rm prescribed}\quad
H\in{\mathbf H}_n~.
\end{equation}
Adopting the convention that $\overline{\mathbf P}_n$ is the {\em closure} (in
$\overline{\mathbf H}_n$) of the {\em open} set ${\mathbf P}_n~$, one can
say that $\overline{\mathbf L}_H$ is the {\em closure} of the {\em open} set
$\mathbf{L}_H~$.
\vskip 0.2cm

\noindent
The following is well known, see e.g. \cite[Section 5.1]{BGFB1994}
and for a special case \cite{CohenLewRod1997}.

\begin{Ob}\label{Ob:ExpStab}
If for some $-H\in\mathbf{P}_n$ one has that $\mathbf{L}_H$ from Eq.
\eqref{eq:DefL_H}, is so that 
\begin{equation}\label{eq:QuadStab}
\mathbf{A}\subset\mathbf{L}_H~,
\end{equation}
then one can find ${\scriptstyle\alpha}>0$ and ${\scriptstyle\beta}\geq 1$ so
that the solution $x(t)$ of the equation in \eqref{eq:DiffInclusion}
uniformly satisfies,
\begin{equation}\label{eq:ExpStab}
{\scriptstyle\beta}\|x(t_o)\|e^{\alpha(t_o-t)}\geq\|x(t)\|
\quad\quad\forall x(t_o)\quad\forall t\geq t_o\geq 0.
\end{equation}
\end{Ob}

\noindent
The celebrated ~{\em Linear Matrix Inequality}~ (LMI)
technique\begin{footnote}{where $m$ and $n$ are ``modest".}\end{footnote},
see e.g. \cite{BGFB1994}, \cite{GahNemiLaubChila1995}, is the
prominent engineering tool to finding whether or not for
\mbox{$\mathbf{A}=\left\{A_1,~\ldots~,~ A_m\right\}$} there exists $H$
satisfying \eqref{eq:QuadStab}.
\vskip 0.2cm

\noindent
Already here we need to recall that the converse of Observation \ref{Ob:ExpStab}
is in general not true. Namely Eq. \eqref{eq:ExpStab} does not imply Eq.
\eqref{eq:QuadStab}. For a special case where the two conditions are
equivalent see, \cite{CohenLewRod1997}.
\vskip 0.2cm

\noindent
We next start exploring the structure of the set $\mathbf{L}_H$.

\begin{Dn}\label{Dn:CIC}
{\rm
A set of $n\times n$ matrices is said to be 
a~ {\em Convex Cone}
if it is closed under positive scaling and summation.\\
A set of $n\times n$ matrices is said to be~ {\em Invertible} (=``closed
under inversion") if whenever a matrix $M$ in it is non-singular,
its inverse $M^{-1}$, belongs to the same set.\\
A set of $n\times n$ matrices\begin{footnote}{Convex Invertible Cones were
originally defined over any real unital algebra, see
\cite{CohenLew1997a}, \cite{CohenLew2007}, \cite{LewRodYar2005}. For
simplicity of exposition, we here start with
matrices.}\end{footnote} combining both properties
is called a~ {\em Convex Invertible Cone}, {\bf cic} in short.
\qed
}
\end{Dn}

\begin{Ex}
{\rm 1.
The set $\mathbf{H}_n$ of $n\times n$ {\em non-singular} Hermitian matrices
is a cone, closed under inversion, but not convex as $\pm H$ may belong to
it, but not their sum\begin{footnote}{The set $\overline{\mathbf H}_n$ of
all $n\times n$ Hermitian matrices, will be addressed Observation
\ref{Ob:MatrixConvexity}.}
\end{footnote}.\\
2.
The set of $2\times 2$ matrices with ${\rm det}=1$ is closed under inversion,
but not convex. Its convex subset of matrices of the form
$\left(\begin{smallmatrix} 0&-1\\1&~~c\end{smallmatrix}\right)$,
${\scriptstyle c}\in\C$, is not closed under inversion, as
\mbox{
$\left(\begin{smallmatrix} 0&-1\\1&~~c\end{smallmatrix}\right)^{-1}=
\left(\begin{smallmatrix}~~c&~1\\-1&~0\end{smallmatrix}\right)$.}\\
3. The set $\overline{\mathbf P}_n$ is
a convex invertible cone, although it contains singular matrices.
\qed
}
\end{Ex}
\vskip 0.2cm

\noindent
The first fundamental structural result is the
following, see e.g. \cite[Lemma 3.5, Proposition 3.7]{CohenLew1997a}.

\begin{Tm}\label{Tm:L_H=CIC}
For arbitrary $H\in\mathbf{H}_n$, the set $\mathbf{L}_H$ in \eqref{eq:DefL_H} is:
A cone, closed under inversion, contains the matrix $H$ and a maximal open
convex set of non-singular matrices. 
\end{Tm}

\noindent
Specifically, maximality is in the following sense: Whenever
\mbox{$B\not\in\overline{\mathbf L}_H$} (i.e. $B$ does not belong to the
closure of the open set $\mathbf{L}_H$), it means that
\mbox{$\min\limits_{j=1,~\cdots~n}\lambda_j(HB+B^*H)=-\beta$} for some
scalar $\beta>0$, then one can always find $A\in\mathbf{L}_H$ so that $A+B$
is singular.
\vskip 0.2cm

\noindent
Define now a matrix $A:={\scriptstyle\frac{1}{2}}H^{-1}({\beta}I+B^*H-HB)$.
On the one hand $A\in\mathbf{L}_H$ (i.e. $HA\in\mathbf{L}_I$),
on the other hand,
\[
H(A+B)=\underbrace{
{\scriptstyle\frac{1}{2}}({\beta}I+B^*H-HB)}_{HA}+
\underbrace{
{\scriptstyle\frac{1}{2}}(HB+B^*H-B^*H+HB)}_{HB}
=
{\scriptstyle\frac{1}{2}}({\beta}I+HB+B^*H),
\]
which by construction is a singular (Hermitian) matrix. Now as $H$ is
non-singular, \mbox{$H(A+B)$} is singular if and only if, 
\mbox{$(A+B)$} is singular.
\vskip 0.2cm

\noindent
In \cite[Section 3]{Ando2001} T. Ando {\em characterized} the set
$\mathbf{L}_H$ for $H\in\mathbf{P}_n$ and in
\cite[Theorem 3.5]{Ando2004} he extended it to\begin{footnote}{To
be precise, the result is formulated for the
Stein inclusion $\{A\in\C^{n\times n}:~(H-A^*HA)\in\mathbf{P}_n\}$, where
$H\in\mathbf{H}_n$ is prescribed.}\end{footnote} $H\in\mathbf{H}_n$.
In particular, he showed that the conditions in Theorem \ref{Tm:L_H=CIC}
fall short from characterizing the set $\mathbf{L}_H~$, namely the
converse statement is (significantly) more involved.
\vskip 0.2cm

\noindent
Motivated by {\em physical} considerations, in this work we focus
on the special case of the set $\mathbf{L}_H$ where $H=I$.
Consequently, {\em convexity} in Theorem \ref{Tm:L_H=CIC}, can be
substituted by {\em matrix-convexity}, in Theorem \ref{Tm:L_I=P+iH}
below.

\section{matrix-convex sets and cones of matrices}
\label{Sec:Mcic}
\setcounter{equation}{0}

We next resort to the notion of a {\em matrix-convex} set, see e.g.
\cite{EffrWink1997} and more recently, \cite{EverHeltKlepMcCull2018},
\cite{Kriel2019} and \cite{PassShalSol2018}.

\begin{Dn}\label{Dn:MatrixConvex}
{\rm
{\bf a.}~ A family of square matrices\begin{footnote}{We do not assume
that $\mathbf{A}\subset\overline{\mathbf H}$.}\end{footnote} (of various
dimensions) $\mathbf{A}$, is said to be~ {\em matrix-convex of level
$n$}, if for all $\nu=1,~\ldots~,~n$,\\
for all natural $k$, 
\begin{equation}\label{eq:Isometry1}
\sum\limits_{j=1}^k\upsilon_j^*\upsilon_j=I_{\nu}
\quad\quad
\begin{smallmatrix}
\forall{\upsilon}_j\in\C^{{\eta}_j\times\nu}
\\~\\
\eta_j\in[1,~\nu],
\end{smallmatrix}
\end{equation}
one has that having $A_1,~\ldots~,~A_k$ (of various dimensions
\mbox{$1\times 1$} through \mbox{$\nu\times\nu$}) within
$\mathbf{A}$, implies that
\begin{equation}\label{eq:DefMatrixConvex}
\sum\limits_{j=1}^k\upsilon_j^*A_j\upsilon_j~,
\end{equation}
belongs to $\mathbf{A}$ as well.
\vskip 0.2cm

\noindent
If the above holds for all $n$, we say that the set $\mathbf{A}$ is
{\em matrix-convex}.
\vskip 0.2cm

\noindent
{\bf b.} ~
A family of square matrices $\mathbf{A}$ is said to be~ {\em matrix-convex
cone} if the right hand side of Eqs. \eqref{eq:Isometry1}
is relaxed to be in $\mathbf{P}_{\nu}$, with $\nu=1, 2,~\ldots$
\qed
}
\end{Dn}
\vskip 0.2cm

\noindent
Matrix-convex cones are closely related to the classical notion of~ {\em
Complete Positivity},~ see e.g. \cite{Ando1991}, \cite{Choi1975}, and
for a comprehensive account of the subject, see \cite{BeranmShaked2003}.
In recent years it has been applied to the study of {\em Quantum
Channels}, see e.g. \cite{Levick2018}.
\vskip 0.2cm

\noindent
We next present some prime examples of matrix-convex sets and cones. To this
end, recall that we denote by $\overline{\mathbf H}$  the set of (possibly
singular) Hermitian matrices. Skew-Hermitian matrices are denoted by,
$i\overline{\mathbf H}$. It is common to consider $\overline{\mathbf H}$
and $i\overline{\mathbf H}$ as the matricial extensions of $\R$ and $i\R$,
respectively.

\begin{Ob}\label{Ob:MatrixConvexity}
(I)\quad 
Each of the following families of matrices,
\[
\overline{\mathbf H}~,\quad\quad\quad
i\overline{\mathbf H}~,\quad\quad\quad
\overline{\mathbf P},\quad\quad\quad
{\mathbf P},
\]
is a matrix-convex cone.
\vskip 0.2cm

\noindent
(II)\quad 
$\overline{\mathbf A}_2({\scriptstyle\alpha})$,
(${\mathbf A}_2({\scriptstyle\alpha})$), the closed (open) family of
square matrices whose spectral norm is uniformly bounded (with a
prescribed ${\scriptstyle\alpha}>0$),
\begin{equation}\label{eq:A_2}
\begin{matrix}
\overline{\mathbf A}_2({\scriptstyle\alpha})&=&
\{~A~:~{\scriptstyle\alpha}\geq\| A\|_2~\},
\\~\\
{\mathbf A}_2({\scriptstyle\alpha})&=&
\{~A~:~{\scriptstyle\alpha}>\| A\|_2~\},
\end{matrix}
\end{equation}
is a matrix-convex set. 
\end{Ob}

\noindent
Verification of (I) and (II) is self-evident and thus omitted.
\vskip 0.2cm

\noindent
In Observation \ref{Ob:MatrixConvexity} the sets $\overline{\mathbf H}_n$
or $\mathbf{P}_n$ are replaced by $\overline{\mathbf H}$ or $\mathbf{P}$
indicating that matrix-convexity is in principle dimension-free.
\vskip 0.2cm

\noindent
Substituting in Eq. \eqref{eq:A_2} ${ \alpha}=1$, one obtains the
matrix-convex (closed) open contractions ($\overline{\mathbf A}_2(1)$)
and ${\mathbf A}_2(1)$, which are pivotal to {\em discrete-time}
passivity, see \cite{Lewk2020c}.
\vskip 0.2cm

\noindent
\begin{Rk}\label{Rk:MatrixConvex}
{\rm
(i)\quad Substituting in Eq. \eqref{eq:DefMatrixConvex} $k=1$, reveals
that matrix-convexity in particular implies that the set $\mathbf{A}$
is invariant under unitary similarity.
\vskip 0.2cm

\noindent
Thus in particular, in Eq. \eqref{eq:A_2} the spectral norm $\|~\|_2$,
can not be substitute d by another unitarily-variant, induced
matrix norm e.g.  $\|~\|_1$ or $\|~\|_{\infty}~$.
\vskip 0.2cm

\noindent
(ii)\quad Similarly, taking in Eq. \eqref{eq:DefMatrixConvex},
\mbox{${\upsilon}_j=r_j{I}_{{\nu}_j}$} where \mbox{$r_1~,~\ldots~,~r_k$}
are arbitrary real scalars so that
${r_1}^2+~\ldots~+{r_k}^2=1$,
reveals that matrix-convexity in particular implies classical convexity.
\vskip 0.2cm

\noindent
(iii)\quad In light of the two above items, we here show convexity
combined with closure under unitary similarity, fall short from
implying matrix-convexity: Consider the set of positively scaled
identity matrices, ${\scriptstyle \alpha}I~:~{\scriptstyle\alpha}>0~\}$.
Trivially, this set is convex and each matrix is invariant under unitary
similarity. Yet this set is not {\em matrix-convex}. Indeed, already
for $k=2$ and arbitrary
$~{\scriptstyle\alpha}_1\not={\scriptstyle\alpha}_2$,
\[
\left(\begin{smallmatrix}1&0\\~\\0&0_{1\times({\nu}_1-1)}\end{smallmatrix}\right)
{\alpha}_1I_{{\nu}_1}
\left(\begin{smallmatrix}1&0\\~\\0&0_{1\times({\nu}_1-1)}\end{smallmatrix}\right)
+
\left(\begin{smallmatrix}0_{1\times({\nu}_2-1)}&0\\~\\0&1\end{smallmatrix}\right)
{\alpha}_2I_{{\nu}_2}
\left(\begin{smallmatrix}0_{1\times({\nu}_2-1)}&0\\~\\0&1\end{smallmatrix}\right)
=
\left(\begin{smallmatrix}{\alpha}_1&&0\\~\\0&&{\alpha}_2\end{smallmatrix}\right).
\]
Hence, the set of scaled identity matrices is not matrix-convex.
\qed
}
\end{Rk}

matrix-convexity

\section{Maximal non-singular matrix-convex cones, closed under inversion}
\label{Sec:Maximal Mcic}
\setcounter{equation}{0}

\noindent
As already mentioned, we here focus on the special case of the sets
$\mathbf{L}_H$ ($~\overline{\mathbf L}_H$) in Eq. \eqref{eq:DefL_H},
where one substitutes $H=I$, i.e.
\begin{equation}\label{eq:DefL_I}
\begin{matrix}
\mathbf{L}_I&:=&\{ A~:~A+A^*\in\mathbf{P}~\}
\\~\\
\overline{\mathbf L}_I&=&\{ A~:~A+A^*\in\overline{\mathbf P}~\}.
\end{matrix}
\end{equation}
Note now that the sets in Eq. \eqref{eq:DefL_I} may be viewed as a
matricial extensions of $\C_R$, $\overline{\C}_R$, respectively. Indeed,
one can equivalently write these sets as,
\begin{equation}\label{eq:AlternativeL_I}
\begin{matrix}
\mathbf{L}_I&=&\{P+iH~:~P\in\mathbf{P}~,~H\in\overline{\mathbf H}\}
\\~\\
\overline{\mathbf L}_I&=&
\{P+iH~:~P\in\overline{\mathbf P}~,~H\in\overline{\mathbf H}\}.
\end{matrix}
\end{equation}
Under the assumption $H=I$, the stability of differential inclusion in
Observation \ref{Ob:ExpStab} takes the form that Eq. \eqref{eq:ExpStab}
holds with $\beta=1$ and the norm used is the spectral norm
(i.e. $\| x\|_2=\sqrt{x^*x}$).
\vskip 0.2cm

\noindent
Here is the first motivation to resorting to the notion of
matrix-convexity.

\begin{Tm}\label{Tm:L_I=P+iH}
The following statements are true.

\begin{itemize}
\item[(i)~~~]{}
The set ${\mathbf L}_I$ in Eq. \eqref{eq:DefL_I} is: A cone, closed
under inversion, contains the matrix $I$ and a maximal open
matrix-convex set of non-singular matrices.
\vskip 0.2cm

\item[(ii)~~]{}
Conversely, a cone closed under inversion, containing the matrix $I$
and a maximal open convex set of non-singular matrices, which in
addition is closed under unitary similarity. is the set ${\mathbf L}_I$.
\vskip 0.2cm

\item[(iii)~]{}
The set $\overline{\mathbf L}_I$ in \eqref{eq:DefL_I} is a cone closed
under inversion and a closed matrix-convex
set containing the matrix $I$, and on its boundary the matrix $~iI$.
\vskip 0.2cm

\item[(iv)~]{}
$\overline{\mathbf L}_I\bigcap\overline{\mathbf L}_{-I}=
i\overline{\mathbf H}$. The set $~i\overline{\mathbf H}~$ is a
matrix-convex cone, closed under inversion; in fact, a maximal convex
subset of $~\C^{n\times n}$, which does not contain an involution.
\end{itemize}
\end{Tm}

\noindent
{\bf Proof :}~(i)\quad Almost all properties are obtained from Theorem
\ref{Tm:L_H=CIC}, upon substituting $H=I$. It is only matrix-convexity
which is left to be shown.
\vskip 0.2cm

\noindent
Recall that from Eq. \eqref{eq:AlternativeL_I} it follows that for
$j=1,~\ldots~,~k$ having $A_j\in\mathbf{L}_I$ means that
$A_j=P_j+iH_j$
for some $P_j\in\mathbf{P}$ and some $H_j\in\overline{\mathbf H}$.
In principle, matrix-convexity can now be deduced from from item (I) of
Observation \ref{Ob:MatrixConvexity}. Still, we here write it explicitly,
for $~A_j\in\mathbf{L}_I$,
\[
\sum\limits_{j=1}^k\upsilon_j^*A_j\upsilon_j=
\sum\limits_{j=1}^k\upsilon_j^*\underbrace{(P_j+iH_j)}_{A_j}\upsilon_j=
(~
\underbrace{(\sum\limits_{j=1}^k\upsilon_j^*P_j\upsilon_j)}_{\in\mathbf{P}}+
i(\sum\limits_{j=1}^k\upsilon_j^*H_j\upsilon_j)~)\in\mathbf{L}_I,
\]
so this part of the claim is established.
\vskip 0.2cm

\noindent
(ii)\quad 
We shall show the converse in steps.
\vskip 0.2cm

\noindent
Step 1:\quad
We first show that the spectrum of each matrix in the set,
is in $\C_R~$.\\
Assume that there exists a non-singular matrix $A$ in this set,
so that $(A-{\lambda}I)v=0$ for some\begin{footnote}{As before,
matrix-convexity formulation is dimension-free. The introduction of
$0\not=v\in\C^n$ is only to simplify the proof presentation.}\end{footnote}
$0\not=v\in\C^n$ and some $\lambda\in\overline{\C}_L~$. Then let
$A_1$ be the following convex combination of $A$ and $A^{-1}$,
\[
A_1={\scriptstyle\frac{1}{1+|\lambda|^2}}A+
{\scriptstyle\frac{|\lambda|^2}{1+|\lambda|^2}}A^{-1}.
\]
By construction $(A_1-2{\rm Re}(\lambda)I)v=0$ with the same
$v$.
\vskip 0.2cm

\noindent
Now if $\lambda$ was on $i\R$, then ${\rm Re}(\lambda)=0$ so $A_1$
is singular and we are done. Hence, assume that $\lambda$ was in $\C_L$,
i.e. $0>{\rm Re}(\lambda)$. Then for 
\mbox{$\alpha=-2{\rm Re}(\lambda)>0$}, with the same $v$,
\mbox{$(A_1+{\alpha}I)v=0$}, so the matrix \mbox{$A_1+{\alpha}I$}
is singular.
\vskip 0.2cm

\noindent
Assume from now on that the spectrum of each matrix in the set,
is in $\C_R~$.
\vskip 0.2cm

\noindent
Step 2:\quad
Here are some facts on the sets
of the form $\mathbf{L}_H$ in Eq. \eqref{eq:DefL_H} where
$H\in\mathbf{P}_n~$.\\
(a)~ $\bigcup\limits_{H\in\mathbf{P}}\mathbf{L}_H$, covers all matrices
whose spectrum is in $\C_R~$. In particular, each matrix whose spectrum
is in $\C_R$ belongs to (infinitely) many sets of the form
$\mathbf{L}_H$.\\
(b)~ Each set of the form $\mathbf{L}_H$ with
$H\in\mathbf{P}_n$ contains the matrix $I$.\\
(c)~ By Theorem \ref{Tm:L_H=CIC}, each of the sets $\mathbf{L}_H$,
$H\in\mathbf{P}$, is
a maximal open convex cone of matrices whose spectrum is in $\C_R~$.\\
(d)~ All these sets are similar, in the sense that,
\mbox{$H^{\frac{1}{2}}\mathbf{L}_{H}H^{-\frac{1}{2}}
=\mathbf{L}_I~$,} see e.g.  \cite[Lemma 3.4]{CohenLew1997a}.
\vskip 0.2cm

\noindent
Step 3:\quad
It is left to check which of the sets of the form $\mathbf{L}_H$ with
$H\in\mathbf{P}_n~$, is closed under unitary similarity.\\
Recall (e.g. \cite[Lemma 3.4]{CohenLew1997a}) that for a unitary matrix
$U$ and $H\in\mathbf{H}_n$, both arbitrary, one has that
\mbox{$U^*\mathbf{L}_HU=\mathbf{L}_{U^*HU}$}. Thus to guarantee that
\mbox{$U^*\mathbf{L}_HU=\mathbf{L}_H$}, one must take $H=\pm{I_n}$.
To conform with previous steps, $H=I$, so this item is established.

\noindent
(iii)\quad This follows from the previous items along with the already mentioned
fact that $\overline{\mathbf L}_H$ is the closure of the open set $\mathbf{L}_H$
for all $H\in\mathbf{H}$, and in particular for $H=I$.
\vskip 0.2cm

\noindent
(iv)\quad 
See item $(I)$ of Observation \ref{Ob:MatrixConvexity} along with
\cite[Proposition 3.2.5(i)]{CohenLew2007}.
\qed
\vskip 0.2cm

\noindent
We conclude this section by putting Theorem
\ref{Tm:L_I=P+iH} into perspective:
\vskip 0.2cm

\noindent
{\bf 1.}~~In
\cite[Theorem 3.4]{AlpayLew2020}, and further in \cite[Section5]{Lewk2020d},
a quantitative refinement of items (i), (ii) of Theorem \ref{Tm:L_I=P+iH},
is introduced, where the (linear) {\em Lyapunov} inclusion in Eq.
\eqref{eq:DefL_H} is substituted by a (quadratic) {\em Hyper-Lyapunov}
inclusion.
\vskip 0.2cm

\noindent
{\bf 2.}~~As already mentioned, a complete characterization of the set
$\mathbf{L}_H$ for an arbitrary $H\in\mathbf{H}_n$, appeared in
\cite[Theorem 3.5]{Ando2004}. The restriction in Theorem
\ref{Tm:L_I=P+iH} to $H=I$ enables us, by resorting to the notion of
matrix-convexity, to obtain a much simpler description of $\mathbf{L}_I$,
which in turn is exploited in presenting Positive Real functions, see
Definition \ref{Dn:PosReal} below.

\section{Maximal matrix-convex invertible cones of Rational Functions}
\label{MCICs of Rational Functions}
\setcounter{equation}{0}

\noindent
As a first connection with the structure we focus on, we cite the following
adapted version of \cite[Proposition 5.3.2]{CohenLew2007}.

\begin{Pn}\label{Pn:cic(f,g)}
Let $f(s)=\frac{1}{s}$ and $g(s)\equiv 1$ be a pair of scalar rational positive
real functions, of degree 1 and 0, respectively. 
\vskip 0.2cm

\noindent
A scalar positive real rational function can always be generated
by iteratively taking positive scaling, summation and inversion 
of $f(s)$ and $g(s)$.
\end{Pn}

\noindent
Thus, one can conclude that scalar rational $\mathcal{PR}$ functions may be viewed
as ${\rm\bf cic}\left(f, g\right)$ a convex invertible cone generated by the above
$f(s)$ and $g(s)$. An analogous observation for state-space realization of the
above $f$ and $g$, will be given in Example \ref{Ex:RealizationArray} below.
\vskip 0.2cm

\noindent
The fact that in the scalar case, matrix-convexity degenerates to classical
convexity, simplified the above discussion. We now proceed to
matrix-valued rational functions.
\vskip 0.2cm

\noindent
Recall that $\overline{\mathbf L}_I$ is the matricial generalizations of 
$\overline{\C}_R$. Thus, we find it convenient to employ the terminology of Eq.
\eqref{eq:DefL_I} to describe matrix-valued Positive Real functions.

\begin{Dn}\label{Dn:PosReal}
{\rm 
Let $F(s)$ be an $m\times m$-valued rational function so that
\mbox{${F(s)}_{|_{s\in\R}}\in\mathbb{R}^{m\times m}$}.  $F(s)$ is said to
be~ {\em Positive Real}, denoted by $\mathcal{PR}$, if it analytically maps
$~\C_R$ to $~\overline{\mathbf L}_{I_m}\bigcup\infty$.
}
\qed
\end{Dn}

\noindent
As already mentioned, this set of $m\times m$-valued positive real
rational functions corresponds to the driving point immittance of a lumped
\mbox{$R-L-C$} electrical networks, (along with transformers) with $~m~$
inputs and $~m~$ outputs. To be precise, in the non-reciprocal case,
gyrators are needed as well.
\vskip 0.2cm

\noindent
Here is a fundamental structural property of this set.

\begin{Tm}\label{Tm:MaximalMatrixPR}
The family $\mathcal{PR}$, of $m\times m$-valued positive real rational
functions, is a cone, closed under inversion and a maximal matrix-convex 
family of functions which is analytic in $\C_R$.
\vskip 0.2cm

\noindent
Conversely, a maximal matrix-convex cone of $m\times m$-valued rational
functions, analytic in $\C_R$, containing the zero degree function
\mbox{$F_o(s)\equiv I_m$}, is the set $\mathcal{PR}$.
\end{Tm}

\noindent
{\bf Proof :}\quad
Using the fact that all functions in $\mathcal{PR}$ analytically map $\C_R$ to
$\overline{\mathbf L}_{I_m}$, together with item (i) of Theorem
\ref{Tm:L_I=P+iH} establishes the sought structure.
\vskip 0.2cm

\noindent
For maximality, take a rational function $G(s)$ which does not belong to
$\mathcal{PR}$. To avoid triviality, assume that it is analytic in $\C_R$,
but there exists $s_o\in\C_R$ so that
\mbox{$\left(G(s_o)\right)v=(-a+i\beta)v$}
for some $a>0$, $\beta\in\R$ and $0\not=v\in\C^m$.
\vskip 0.2cm

\noindent
Note now that (with the same $v$),
$\underbrace{\left(G(s)+aI_m\right)_{|_{s_o}}}_{M\in\C^{m\times m}}v=i\beta{v}$
and thus, $(M+{\beta}^2M^{-1})v=0$, i.e.
\[
(~\overbrace{G(s)+aI_m}^M+{\beta}^2\overbrace{
\left(G(s)+aI_m\right)^{-1}}^{M^{-1}}~)_{|_{s_o}}
v=0.
\]
This means that the rational function
\mbox{$G(s)+aF_o(s)+{\beta}^2\left(G(s)+aF_o(s)\right)^{-1}$},
is within ${\rm\bf cic}(G, F_o)$ (see Definition \ref{Dn:CIC}), and
has a right-half plane zero at $s_o$ and thus its inverse,
\[
\left(G(s)+aF_o(s)+{\beta}^2\left(G(s)+aF_o(s)\right)^{-1}\right)^{-1},
\]
(still within ${\rm\bf cic}(G, F_o)$) is no longer analytic in $\C_R$.
Hence, the claim is established.
\qed
\vskip 0.2cm

\noindent
Note that in principle the same proof applies to not necessarily rational
$\mathcal{PR}$ function. A scalar version of Theorem
\ref{Tm:MaximalMatrixPR} appeared in \cite[Proposition 4.1.1]{CohenLew2007}.
\vskip 0.2cm

\noindent
We here point out that in \cite{AlpayLew2020} a subset of $\mathcal{PR}$ 
functions which are associated with absolute stability (a.k.a the Lurie
problem) is studied. In scalar terminology, these are functions mapping
$\C_R$ into a bounded disk within $\C_R$. Furthermore, under inversion
this disk is mapped onto itself.
\vskip 0.2cm

\noindent
We conclude this section by illustrating an application of Theorem
\ref{Tm:MaximalMatrixPR}. Here are the details. 
\vskip 0.2cm

\begin{figure}[ht!]
\begin{tikzpicture}[scale=1.8]
    \draw[color=black, thick]
       (1,0) to [short,o-] (4,0){} % Baseline for connection to ground
       (0.9,1) node[]{\large{$\mathbf{Z_{\rm in}~~\rightarrow}$}}
       (1,2) to [short,o-] (2,2)
       (3,2)   to node[short]{} (4,2)
       (2,2) to [short] (3,2)
       (2,2) to [C=$C_a$,*-*] (2,1)
       (2,1) to [L=$L_b$, *-*] (2,0)
       (3,0) to [L,l=$L_c$, *-*] (3,2)
       (4,0) to [C, l=$C_d$, *-*] (4,2)
        ;
 \end{tikzpicture}
 \caption{$
{\rm Z}_{\rm in}(s)=\left(\left((sC_a)^{-1}+sL_b\right)^{-1}+
(sL_c)^{-1}+sC_d\right)^{-1}.$}
\label{Fig:Circuit2FeedbackLoops}
\end{figure}
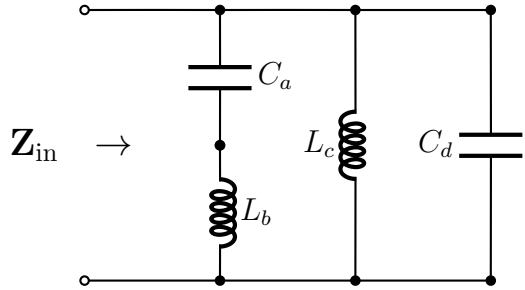

\noindent
The driving point impedance of the circuit in Figure
\ref{Fig:Circuit2FeedbackLoops}, is a standard positive real (odd a.k.a.
lossless or Foster) rational function of degree four. Yet,
employing the notation,
\begin{equation}\label{eq:Phi}
\phi(X, Y):=\left(X^{-1}+Y\right)^{-1},
\end{equation}
this driving point impedance can also be written as,
\begin{equation}\label{eq:EquivZin}
Z_{\rm in}(s)=\phi\left(~\phi(sC_a, sL_b),~\phi(sL_c,~sC_d)~\right).
\end{equation}
We now leave this circuit for a short while and address a $2m\times 2m$-valued
feedback-loop network $H(s)$ in Figure \ref{Fig:MultipleFeedbackLoop}.
\begin{figure}[ht!]
\begin{tikzpicture}[scale=1.1]
\draw[color=black, thick] [->] (1.5,0.5) -- (3,0.5){};
\draw[color=black, thick] [->] (4,0.5) -- (5.5,0.5){};
\draw[color=black, thick] [->] (0,3.5) node[left]{${\rm In}_1$} -- (1.15,3.5) {};
\draw[color=black, thick] [->] (7.0,2.0)node[right]{${\rm In}_2$} -- (5.85,2.0){};
\draw[color=black, thick] [->] (1.85,3.5) -- (3,3.5){};
\draw[color=black, thick] [->] (4,3.5) -- (5.5,3.5){};
\draw[color=black, thick] [->] (5.5,3.5) -- (7.0,3.5) node[right] {${\rm Out}_1$};
\draw[color=black, thick] [->] (1.5,2.0) -- (0.0,2.0) node[left] {${\rm Out}_2$};
\draw[color=black, thick] [->] (5.5,0.5) -- (5.5,1.65){};
\draw[color=black, thick] [->] (1.5,2) -- (1.5,0.5){};
\draw[color=black, thick] [->]  (1.5,2) -- (1.5,3.15){};
\draw[color=black, thick] [->]  (1.5,5) -- (1.5,3.85){};
\draw[color=black, thick] [->]  (3,5) -- (1.5,5){};
\draw[color=black, thick] [->]  (5.5,5) -- (4,5){};
\draw[color=black, thick] [->]  (5.5,3.5) -- (5.5,5){};
\draw[color=black, thick] [->]  (5.5,3.5) -- (5.5,2.35){};
\draw[color=black, thick] [->]  (3,2) -- (1.5,2){};
\draw[color=black, thick] [->]  (5.15,2) -- (4,2){};
\draw[color=black, thick] [-]  (3,5.35) -- (4,5.35){};
\draw[color=black, thick] [-]  (3,5.35) -- (3,4.65){};
\draw[color=black, thick] [-]  (4,5.35) -- (4,4.65){};
\draw[color=black, thick] [-]  (4,5.35) -- (4,4.65){};
\draw[color=black, thick] [-]  (3,4.65) -- (4,4.65){};
\draw[color=black, thick] [-]  (3,3.85) -- (4,3.85){};
\draw[color=black, thick] [-]  (3,3.15) -- (4,3.15){};
\draw[color=black, thick] [-]  (3,3.15) -- (3,3.85){};
\draw[color=black, thick] [-]  (4,3.15) -- (4,3.85){};
\draw[color=black, thick] [-]  (3,2.35) -- (4,2.35){};
\draw[color=black, thick] [-]  (3,2.35) -- (3,1.65){};
\draw[color=black, thick] [-]  (4,2.35) -- (4,1.65){};
\draw[color=black, thick] [-]  (3,1.65) -- (4,1.65){};
\draw[color=black, thick] [-]  (3,0.85) -- (4,0.85){};
\draw[color=black, thick] [-]  (3,0.15) -- (4,0.15){};
\draw[color=black, thick] [-]  (3,0.15) -- (3,0.85){};
\draw[color=black, thick] [-]  (4,0.15) -- (4,0.85){};
\draw[color=black, thick] (1.5,3.5) circle (0.35){};
\draw[color=black, thick] (5.5,2.0) circle (0.35){};
\draw[color=black, thick] (3.5,5)node {$F_d(s)$};
\draw[color=black, thick] (3.5,3.5)node {$F_c(s)$};
\draw[color=black, thick] (3.5,2)node {$F_a(s)$};
\draw[color=black, thick] (3.5,0.5)node {$F_b(s)$};
\draw[color=black, thick] (1.05,3.5)node[below] {+};
\draw[color=black, thick] (1.5,3.95)node[left] {-};
\draw[color=black, thick] (1.5,3.05)node[right] {-};
\draw[color=black, thick] (5.5,2.4)node[left] {+};
\draw[color=black, thick] (5.95,2.0)node[above] {+};
\draw[color=black, thick] (5.5,1.55)node[right] {-};
\end{tikzpicture}
\caption{$2m\times 2m$-valued network}
\label{Fig:MultipleFeedbackLoop}
\end{figure}
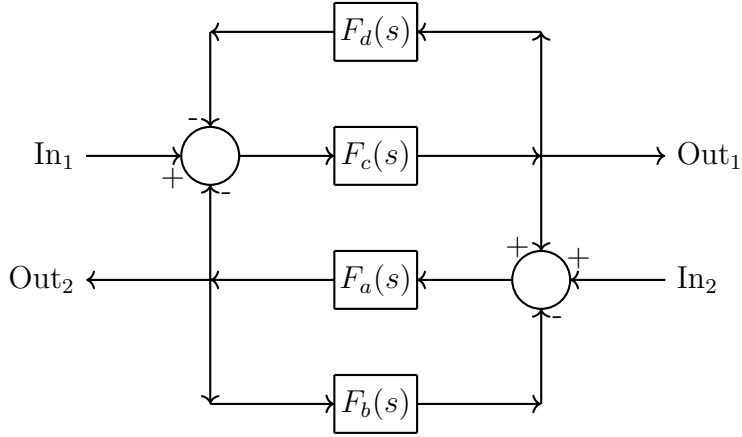
\begin{equation}\label{eq:MultipleFeedbackLoop}
\left(
\begin{smallmatrix}
{\rm Out}_1
\\~\\~\\
{\rm Out}_2
\end{smallmatrix}
\right)
=
\underbrace{
\left(
\begin{smallmatrix}
\left(\hat{F}_c+\hat{F}_a^{-1}\right)^{-1}&-\left(\hat{F}_c+
\hat{F}_a^{-1}\right)^{-1}\hat{F}_a^{-1}
\\~\\
\hat{F}_a^{-1}\left(\hat{F}_c+
\hat{F}_a^{-1}\right)^{-1}&\left(\hat{F}_c^{-1}+\hat{F}_a\right)^{-1}
\end{smallmatrix}
\right)}_{H(s)}
\left(
\begin{smallmatrix}
{\rm In}_1
\\~\\~\\
{\rm In}_2
\end{smallmatrix}
\right)
\quad
\begin{smallmatrix}
\hat{F}_c:=F_c^{-1}+F_d
\\~\\
\hat{F}_a:=F_a^{-1}+F_b~.
\end{smallmatrix}
\end{equation}
On the one hand, $H(s)$ is positive real whenever its four $m\times m$-valued
building blocks $F_a(s)$, $F_b(s)$, $F_c(s)$ and $F_d(s)$ are positive real.
Note however, that the interest in such a feedback-loop network, transcends
the realm of positive real functions. Next note that the upper left corner of
Eq. \eqref{eq:MultipleFeedbackLoop} explicitly says,
\begin{equation}\label{eq:OutOneInOne}
{{\rm Out}_1}_{|_{{\rm In}_2\equiv 0}}=\left(\left(\left(F_a(s)\right)^{-1}+
F_b(s)\right)^{-1}+
\left(F_c(s)\right)^{-1}+F_d(s)\right)^{-1}{\rm In}_1~.
\end{equation}
Employing again the map $\phi$ from Eq. \eqref{eq:Phi}, the relation in
Eq. \eqref{eq:OutOneInOne} can be compactly written as,
\[
{{\rm Out}_1}_{|_{{\rm In}_2\equiv 0}}=
\phi\left(~\phi(F_c, F_d),~\phi(F_a,~F_b)~\right){\rm In}_1~.
\]
Now, in comparison to Eq. \eqref{eq:EquivZin}, one can formally identify the
elements $sC_a$, $sL_b$, $sL_c$, $sC_d$, in Figure
\ref{Fig:Circuit2FeedbackLoops} with the blocks $F_a(s)$,
$F_b(s)$, $F_c(s)$, $F_d(s)$ in Eq. \eqref{eq:OutOneInOne}, respectively. 
\vskip 0.2cm

\noindent
This calls for adapting one of the classical construction schemes of
\mbox{$R-L-C$} circuits, e.g. Brune, Bott-Duffin, Darlington, Foster,
Cauer, etc. see e.g. \cite{AnderVongpa1973}, \cite{BottDuff1949},
\cite{Belev1968},
\cite{Duffin}, \cite{MorelliSmith2019}, \cite{Wohl1969}, to 
introducing a design tool for networks of feedback-loops, more
elaborate than that in Figure \ref{Fig:MultipleFeedbackLoop} (and as
mentioned, the building blocks need not be positive real).
\vskip 0.2cm

\noindent
A word of caution: The passage from one-port circuit design to that
of feedback-loops networks can not be straightforward: Typically
blocks like $F_a(s)$, $F_b(s)$, $F_c(s)$, $F_d(s)$ are
{\em non-commutative}. Hence, one needs to formally introduce positive
real rational functions of say $k$ non-commuting variables, mapping
$\underbrace{\overline{\mathbf L}_{I_n}\times~\cdots~
\times\overline{\mathbf L}_{I_n}}_{k~{\rm times}}$,
to $\overline{\mathbf L}_{I_n}$, where $n$ is a parameter.
Further pursuing this direction is beyond the
scope of this work.

\section{Matrix-convex invertible cones of Realization Arrays}
\label{Sec:Realizations}
\setcounter{equation}{0}

\noindent
The renowned Kalman-Yakubovich-Popov Lemma ties up two representations
of positive real functions: Rational functions and corresponding
state-space realizations.
\vskip 0.2cm

\begin{Tm}\label{Tm:SmallKYP}
Let $F(s)$ be an $m\times m$-valued rational function with no pole at
infinity and let $R_F$ be a corresponding $(n+m)\times(n+m)$ state-space
realization array, i.e.
\begin{equation}\label{eq:Realization}
F(s)=C(sI_n-A)^{-1}B+D
\quad\quad\quad
R_F=
\left({\footnotesize
\begin{array}{c|c}
A&B\\
\hline
C&D
\end{array}}\right).
\end{equation}
If for some $H\in\mathbf{P}_n$ one has that
\begin{equation}\label{eq:OriginalLya}
\left(\begin{smallmatrix}-H&~~0\\~\\0&I_m\end{smallmatrix}\right)
R_F+R_F^*
\left(\begin{smallmatrix}-H&~~0\\~\\0&I_m\end{smallmatrix}\right)
\in\overline{\mathbf P}_{n+m}~,
\end{equation}
then $F(s)$ is positive real.
\vskip 0.2cm

\noindent
If the realization $R_F$ in Eq. \eqref{eq:Realization} is minimal, i.e.
$n$ is the McMillan degree of $F(s)$, then the converse is true as well.
\end{Tm}
\vskip 0.2cm

\noindent
This result first appeared in \cite{AnderMoore1968}. The formulation used
here is due to \cite[Section 5]{Willems2nd1972},
\cite[Section II]{Willems1976}. For further details, see e.g.
\cite{AlpayLew2011},
\cite[Chapters 5, 6]{AnderVongpa1973} and \cite[Subsection 2.7.2]{BGFB1994}.
\vskip 0.2cm

\noindent
The formulation of Theorem \ref{Tm:SmallKYP}, adapted from 
\cite[Section 5]{Willems2nd1972},
employs an elegant idea: To treat the above $(n+m)\times(n+m)$ $R$ as having
two faces\begin{footnote}{Like Janus in the Roman mythology}\end{footnote}:
$(i)$ of an {\em array} and $(ii)$ of a {\em matrix}. This will be further
adopted in Theorem \ref{Tm:n,mMatrixConvexRealizationPR} below.
\vskip 0.2cm

\noindent
The following is classical.
\vskip 0.2cm

\begin{Cy}\label{Cy:BalancedRealization}
In Theorem \ref{Tm:SmallKYP}, up to a change of coordinates
\[
R_F~\longrightarrow~\hat{R}_F:=
\left(\begin{smallmatrix}H^{\frac{1}{2}}&0\\~\\0&I_m\end{smallmatrix}\right)
R_F
\left(\begin{smallmatrix}H^{-\frac{1}{2}}&0\\~\\0&I_m\end{smallmatrix}\right),
\]
Eq. \eqref{eq:OriginalLya} may be substituted by,
\begin{equation}\label{eq:KYPbalanced}
\left(\begin{smallmatrix}-I_n&0\\~\\~0&I_m\end{smallmatrix}\right)
\hat{R}_F+\hat{R}_F^*
\left(\begin{smallmatrix}-I_n&0\\~\\~0&I_m\end{smallmatrix}\right)
=Q,\quad\quad\quad Q\in\overline{\mathbf{P}}_{n+m}~.
\end{equation}
In particular, this is the case when the realization $R_F$ in
Eq. \eqref{eq:Realization} is balanced.
\end{Cy}

\noindent
In \cite[Definition 3]{Willems1976} a balanced realization of positive real
system, satisfying Eq. \eqref{eq:KYPbalanced}, is called ``internally passive".
\vskip 0.2cm

\noindent
Recall also that by definition, balanced realization implies minimality.
However, as before, the passage from Eq. \eqref{eq:OriginalLya} to Eq.
\eqref{eq:KYPbalanced}, does not require minimality of realization.
\vskip 0.2cm

\noindent
To study {\em families} of realization simultaneously satisfying Theorem
\ref{Tm:SmallKYP}, we need to introduce a relaxed version of matrix-convexity
(which is still more strict than classical convexity).

\begin{Dn}\label{Dn:n,mMatrixConvex}
{\rm
(I)~ For all $k$, let $v_j\in\C^{(n+m)\times(n+m)}$,
$j=1,~\ldots~,~k$ be block-diagonal so that
\begin{equation}\label{eq:n,mIsometry}
\sum\limits_{j=1}^k
\underbrace{\left(\begin{smallmatrix}{\upsilon}_{j,n}&0\\~\\0&{\upsilon}_{j,m}
\end{smallmatrix}\right)^*}_{{\upsilon}_j^*}
\underbrace{\left(\begin{smallmatrix}{\upsilon}_{j,n}&0\\~\\0&{\upsilon}_{j,m}
\end{smallmatrix}\right)}_{{\upsilon}_j}
=
\left(\begin{smallmatrix}I_n&&0\\~\\0&&I_m\end{smallmatrix}\right).
\end{equation}
A set $\mathbf{R}$ of $(n+m)\times(n+m)$ matrices is said to be~
$n,m$-{\em matrix-convex}~ if having
${\scriptstyle R_1,~\ldots~,~ R_k}$ in $\mathbf{R}$,
implies that
\[
\sum\limits_{j=1}^k
\underbrace{\left(\begin{smallmatrix}{\upsilon}_{j,n}&0\\~\\0&{\upsilon}_{j,m}
\end{smallmatrix}\right)^*}_{{\upsilon}_j^*}
\underbrace{\left(\begin{smallmatrix}A_j&&B_j\\~\\C_j&&D_j
\end{smallmatrix}\right)}_{R_j}
\underbrace{
\left(\begin{smallmatrix}
{\upsilon}_{j,n}&0
\\~\\
0&{\upsilon}_{j,m}
\end{smallmatrix}\right)}_{{\upsilon}_j},
\]
belongs to ${\mathbf R}$, for all natural $k$ and all block-diagonal
${\upsilon}_j\in\C^{(n+m)\times(n+m)}$.
\vskip 0.2cm

\noindent
(II)~ If Eq. \eqref{eq:n,mIsometry} is relaxed to having
\[
\sum\limits_{j=1}^k
\underbrace{
\left(\begin{smallmatrix}
{\upsilon}_{j,n}&0\\~\\0&{\upsilon}_{j,m}
\end{smallmatrix}
\right)^*}_{{\upsilon}_j^*}
\underbrace{
\left(\begin{smallmatrix}
{\upsilon}_{j,n}&0\\~\\0&{\upsilon}_{j,m}
\end{smallmatrix}
\right)}_{{\upsilon}_j}
=
\left(
\begin{smallmatrix}
P_n&&0\\~\\0&&P_m
\end{smallmatrix}
\right)
\quad\quad\quad
\begin{smallmatrix}
P_n\in\mathbf{P}_n\\~\\
P_m\in\mathbf{P}_m~,
\end{smallmatrix}
\]
then $\mathbf{R}$ is said to be an $n, m$ {\em matrix-convex}.
\qed
}
\end{Dn}
\vskip 0.2cm

\noindent
Strictly speaking to conform with Definition \ref{Dn:MatrixConvex},
the block-diagonal isometries in Eq. \eqref{eq:n,mIsometry} should
involve \mbox{${\upsilon}_{j,\nu}\in\C^{n_j\times\nu}$} and
\mbox{${\upsilon}_{j,\mu}\in\C^{m_j\times\mu}$} where
\mbox{$n_j\in[1,~\ldots~,~\nu]$} and \mbox{$m_j\in[1,~\ldots~,~\mu]$}.
We compromized precision to simplify the notation.
\vskip 0.2cm

\noindent
Here is our first motivation to resorting to the notion of
$n,m$-matrix-convexity,
\vskip 0.2cm

\begin{La}\label{La:LmnMatrixVonvex}
For all $~n, m=1, 2,~\ldots$, the sets 
($\overline{\mathbf L}_{\left(\begin{smallmatrix}
-I_n&0\\0&I_m\end{smallmatrix}\right)}$ and)
$~{\mathbf L}_{\left(\begin{smallmatrix}
-I_n&0\\0&I_m\end{smallmatrix}\right)}$ in \eqref{eq:DefL_H}
are, respectively (closed and) open, $~n,m$-matrix-convex cones,
closed under inversion.
\end{La}
\vskip 0.2cm

\noindent
{\bf Proof :}~
Without loss of generality, a matrix $R$ within
${\mathbf L}_{\left(\begin{smallmatrix}
-I_n&0\\0&I_m\end{smallmatrix}\right)}$ can always be written as
\begin{equation}\label{eq:DetailedR}
\begin{matrix}
R&=&
\left(\begin{smallmatrix}
-P_n+iH_n&-M+Y\\~\\M^*+Y^*&P_m+iH_m+M^*P_n^{-1}M
\end{smallmatrix}\right)
&&
\begin{smallmatrix}
P_n\in{\mathbf P}_n&
P_m\in{\mathbf P}_m\\~\\
H_n\in\overline{\mathbf H}_n&
H_m\in\overline{\mathbf H}_m\\~\\
M, Y\in\mathbb{C}^{n\times m}.&~
\end{smallmatrix}
\end{matrix}
\end{equation}
Now, substituting $R_1~,~\ldots~,~R_k$ of the from of Eq. \eqref{eq:DetailedR}
in Definition \ref{Dn:n,mMatrixConvex} one obtains,
\[
\begin{matrix}
\hat{R}&=&\sum\limits_{j=1}^k\left(\begin{smallmatrix}{\upsilon}_{j,n}&0\\~\\
0&{\upsilon}_{j,m}\end{smallmatrix}\right)^*
\left(\begin{smallmatrix}-P_{j,n}+iH_{j,n}&&-M_j+Y_j\\~\\M_j^*+Y_j^*&&
P_{j,m}+M_j^*P_{j,n}^{-1}M_j+iH_{j,m}\end{smallmatrix}\right)
\left(\begin{smallmatrix}{\upsilon}_{j,n}&0\\~\\0&{\upsilon}_{j,m}
\end{smallmatrix}\right).
\end{matrix}
\]
Next note that,
\[
\left(\begin{smallmatrix}-I_n&0\\~\\0&I_m\end{smallmatrix}\right)\hat{R}
+
\hat{R}^*\left(\begin{smallmatrix}-I_n&0\\~\\0&I_m\end{smallmatrix}\right)
=
2
\underbrace{
\sum\limits_{j=1}^k
\left(\begin{smallmatrix}{\upsilon}_{j,n}&0\\~\\0&{\upsilon}_{j,m}
\end{smallmatrix}\right)^*
\left(\begin{smallmatrix}P_{j,n}&&M_j\\~\\{M_j}^*&&P_{j,m}+M_j^*P_{j,n}^{-1}M_j
\end{smallmatrix}\right)
\left(\begin{smallmatrix}{\upsilon}_{j,n}&0\\~\\0&{\upsilon}_{j,m}
\end{smallmatrix}\right)}_{Q_o},
\]
namely,
\[
\sum\limits_{j=1}^k
\left(\begin{smallmatrix}{\upsilon}_{j,n}&0\\~\\0&{\upsilon}_{j,m}
\end{smallmatrix}\right)^*
\left(\begin{smallmatrix}P_{j,n}&&M_j\\~\\{M_j}^*&&P_{j,m}+M_j^*P_{j,n}^{-1}M_j
\end{smallmatrix}\right)
\left(\begin{smallmatrix}{\upsilon}_{j,n}&0\\~\\0&{\upsilon}_{j,m}
\end{smallmatrix}\right)
=Q_o\in\mathbf{P}_{n+m}
\]
Thus also $\hat{R}$ belongs to $\mathbf{R}$, so the claim is established.
\qed
\vskip 0.2cm

\noindent
We now introduce {\em families of realization arrays}~ associated with rational
functions. Before that, a word of caution: For example, 
$R_1=\left({\footnotesize\begin{array}{c|c}A&B\\ \hline C&D\end{array}}\right)$
and
$R_2=\left({\footnotesize\begin{array}{r|r}A&-B\\ \hline-C&D\end{array}}\right)$
are two realizations of the same rational function. Furthermore, $R_1$
is minimal (balanced) if and only if $R_2$ is minimal (balanced). However,
$R_3={\scriptstyle\frac{1}{2}}(R_1+R_2)=
\left({\footnotesize
\begin{array}{c|c}
A&0\\
\hline
0&D
\end{array}}\right)$ is only a non-minimal realization of a zero degree rational
function $F(s)\equiv D$. Yet, as $(n+m)\times(n+m)$ matrices, if $R_1$ belongs
to $\mathbf{L}_{\left(\begin{smallmatrix}-I_n&0\\0&I_m\end{smallmatrix}\right)}$,
then also $R_2$ and $R_3$ belong to the same set.
\vskip 0.2cm

\noindent
More generally, when considering families of realizations $R$ satisfying
Eq. \eqref{eq:KYPbalanced} as {\em matrices}, one obtains only a proper
subset of
$\overline{\mathbf L}_{\left(\begin{smallmatrix}-I_n&0\\0&I_m
\end{smallmatrix}\right)}$. 
\vskip 0.2cm

\begin{Tm}\label{Tm:n,mMatrixConvexRealizationPR}
Given a family of $m\times m$-valued positive real rational functions of
McMillan degree, of at most, $n$, with no poles at infinity. Consider the
corresponding $(n+m)\times(n+m)$ realizations $R$ in Theorem
\ref{Tm:SmallKYP} and in Eq. \eqref{eq:KYPbalanced}.\\
Then (as matrices), this family of realizations is an $n,m$-matrix-convex
cone, closed under inversion.
\end{Tm}

\noindent
The fact that this is a matrix-convex cone follows from Eqs.
\eqref{eq:Realization}, \eqref{eq:KYPbalanced} along with Lemma
\ref{La:LmnMatrixVonvex}.
\vskip 0.2cm

\noindent
To show that this set is closed under inversion, note that assuming $R$ is
non-singular, multiplying Eq. \eqref{eq:KYPbalanced} by $(R^*)^{-1}$ and
$R^{-1}$ from the left and from the right respectively, the resulting
right-hand side is
\mbox{$\hat{Q}:=(R^*)^{-1}QR^{-1}$}. Now if $Q$ is in
$\overline{\mathbf P}_{n+m}$, then so is $\hat{Q}$.
\vskip 0.2cm

\noindent
Note that as before, in Theorem \ref{Tm:n,mMatrixConvexRealizationPR}
we have not assumed minimality of realizations.
\vskip 0.2cm

\noindent
A different formulation (and a different proof) of
this result, appeared in \cite[Proposition 5.3]{Lewk2020b}.
\vskip 0.2cm

\noindent
To illustrate an application of Theorem \ref{Tm:n,mMatrixConvexRealizationPR}
we next show how a set of of realization arrays, may be parametrized by
a pair of representatives.
\vskip 0.2cm

\noindent
\begin{Ex}\label{Ex:RealizationArray}
{\rm
Recall that in Proposition \ref{Pn:cic(f,g)} we stated that scalar positive
real rational functions can be equivalently described as ${\bf cic}(f, g)$,
with $f(s)=\frac{1}{s}$ and $g(s)\equiv 1$.\quad Let now,
\[
R_f=\left({\footnotesize\begin{array}{c|c}
0&1\\ \hline 1&0 \end{array}}\right)
\quad\quad{\rm and}\quad\quad
R_g=\left({\footnotesize\begin{array}{c|c}
0&0\\ \hline 0&1 \end{array}}\right),
\]
be their (balanced) realizations, respectively.
\vskip 0.2cm

\noindent
Treating $R_f$, $R_g$ as matrices, taking positive scaling, summation and
inversion, one obtains $R_h$,
\[
R_h=dR_g+{\beta}^2\left(aR_g+{\beta}R_f\right)^{-1}=
\left({\footnotesize
\begin{array}{c|c}
-a&\beta\\
\hline
\beta&d
\end{array}}\right)
\quad\quad\quad a, \beta, d\geq 0,
\]
which in turn is a (balanced) realization of the function $h(s)$,
\[
h(s)=\frac{{\beta}^2}{s+a}+d
\quad\quad\quad a, \beta, d\geq 0.
\]
Recall now that in Eq. \eqref{eq:Deg1} we pointed out that $h(s)$ is a
parametrization of all positive real rational functions of degree of at most
one, with no pole at infinity.
}
\qed
\end{Ex}
\vskip 0.2cm

\noindent
We conclude by summarizing in brief, that together with \cite{Lewk2020c},
we have shown above that,
\vskip 0.2cm

\begin{center}
\begin{tabular}{|c|c|}
\hline
\multicolumn{2}{|c|}
{\rm Passive linear time-invariant systems and maximal matrix-convexity}
\\
\hline
{\rm discrete-time}&{\rm continuous-time}\\
\hline
{\rm a set closed under product among its elements}&
{\rm a cone closed under inversion}\\
\hline
\end{tabular}
\end{center}
\vskip 0.2cm

\begin{center}
Acknowledgement
\end{center}

\noindent
The fact that the review was thorough and constructive,
is very well appreciated.

\end{document}